\documentstyle[12pt]{article}
\begin{document}

\baselineskip6mm
\newcommand{\ot}{\otimes}
\newcommand{\x}{{\underline x}}
\newcommand{\omfett}{\mbox{\boldmath $\omega$}}
\newcommand{\efett}{\mbox{\boldmath $e$}}
\newcommand{\gfett}{\mbox{\boldmath $g$}}
\newcommand{\xf}{\mbox{\boldmath $x$}}
\newcommand{\alfett}{\mbox{\boldmath $\alpha$}}
\newcommand{\de}{\delta}
\newcommand{\ga}{\gamma}
\newcommand{\e}{\epsilon}
\newcommand{\th}{\theta}
\newcommand{\tmod}{{\cal T}}\newcommand{\amod}{{\cal A}}
\newcommand{\bemod}{{\cal B}}\newcommand{\cmod}{{\cal C}}
\newcommand{\dmod}{{\cal D}}\newcommand{\hmod}{{\cal H}}
\newcommand{\s}{\scriptstyle}
\newcommand{\tr}{{\rm tr}}
\newcommand{\einsop}{{\bf 1}}
\newcommand{\tnull}{\tmod_{\s 0}}
\newcommand{\be}{\begin{equation}\label}\newcommand{\ee}{\end{equation}}
\newcommand{\ba}{\begin{eqnarray*}}\newcommand{\ea}{\end{eqnarray*}}
\def\eqa#1#2{\\\parbox{14.5cm}{\ba #2 \ea}\hfill\parbox{1.5cm}{\be{#1}\ee}}
\newcommand{\ban}{\begin{eqnarray}\label} \newcommand{\ean}{\end{eqnarray}}
\newcommand{\arb}{\begin{array}} \newcommand{\are}{\end{array}}

\newcommand{\iv}{^{-1}}
\newcommand{\ul}{\underline}
\newcommand{\ua}{{\underline\alpha}}
\newcommand{\ub}{{\underline\beta}}
\newcommand{\uu}{{\underline u}}
\renewcommand{\theequation}{\mbox{\arabic{section}.\arabic{equation}}}
\newcommand{\la}{\langle\,}\newcommand{\ra}{\,\rangle}
\newcommand{\mi}{\,|\,}

\def\einsop{{\bf 1}}

\title{ Difference Equations and Highest Weight Modules of $U_q[sl(n)]$  }
\author{
A. Zapletal$^{1,2}$ 
\\{\small\it Institut f\"ur Theoretische Physik}\\
{\small\it Freie Universit\"at Berlin, Arnimallee 14, 14195 Berlin, Germany} 
}
\date{\small\today}

\maketitle
\footnotetext[1]{email: zapletal@physik.fu-berlin.de}
\footnotetext[2]{supported by DFG: Sonderforschungsbereich 288 
     'Differentialgeometrie und Quantenphysik'}
\begin{abstract}
The quantized version of a discrete Knizhnik-Zamolodchikov system is solved 
by an extension of the generalized Bethe Ansatz.  
The solutions are constructed to be of highest weight which means they 
fully reflect the internal quantum group symmetry. 
  
\end{abstract}

\section{Introduction}\label{s1}
This article can be considered as an addendum to the article  
\cite{bkz1} on matrix difference equations and a generalized version of the  
Bethe ansatz. For an introduction to their r\^{o}le in mathematical physics 
the reader is referred to \cite{bkz1}-\cite{FR} and references contained 
therein. 

Though q-deformations of discrete Knizhnik-Zamolodchikov equations have been 
treat\-ed in much detail over the last years \cite{FR}\cite{resh}
it became not completely clear how that solutions are related to the 
underlying symmetry of such problems.

The conventional algebraic formulation of the Bethe ansatz 
demonstrates the close relation between the eigenvector problem 
and the representation theory of its connected symmetry group (either 
classical or q-deformed): Bethe vectors can be constructed to be highest 
weight vectors of irreducible representations and therefore by simply 
counting them one makes certain on spanning the whole space of states. 

However one has to be careful when moving from classical Lie algebras to a 
quantum group, as it can be seen when an 1 dimensional periodic XXX-Heisenberg 
chain is deformed to the anisotropic XXZ-model. 
Deforming the Hamiltonian in a straightforward way won't preserve the 
(quantum) symmetry. Instead one is forced to change the boundary conditions 
\cite{nv} or to take additional terms (arising from the nontrivial toroidal 
topology) into account as done in \cite{kz1}.

The behavior of the difference equation
\be{1.1}
Q(\x;i)\,f(x_1,\ldots,x_i,\ldots,x_N)=f(x_1,\ldots,x_i+\kappa,\ldots,x_N),
\quad i=1,\ldots,N;
\ee
where $f(\x)$ is a vector-valued function on $N$ variables $x_i$, $Q(\x;i)$ 
a family of linear operators and $\kappa$ an arbitrary shift parameter,
indeed resembles this problematic nature: The operators $Q(\x;i)$ can be 
regarded as a sort of generalized transfer matrices and therefore the analogy 
to a quantum spin chain becomes obvious. 

In Section 2 we formulate this equation in a way adapted to quantum symmetry
and obtain solutions by a generalized Bethe ansatz. In Section 3 they are 
shown to be highest weight vectors, additionally we calculate their weights. 
For sake of transparency both sections are fixed to $U_q[sl(2)]$  containing 
all essential features of a quantum group. Finally, for completeness,  
we briefly comprehend  the aspects of the higher ranked case in Section 4, 
followed by a summary of the given results. 
\section{The generalized Bethe ansatz}
\setcounter{equation}{0}
Consider N vector spaces $V_i \simeq \mbox{${\bf C}^2$} $, each given as the 
representation space of the fundamental representation of $U_q[sl(2)]$. The 
basis vectors will be denoted by $|1\rangle$ resp.~$|2\rangle$.
The R-matrix then acts as a linear operator on two of such spaces 
$V_i$ and $V_j$:
\be{2.1}
R_{ij}:\,V_i \ot V_j \to V_j \ot V_i,
\ee
and is given by the quasitriangular Hopf algebraic structure of $U_q[sl(2)]$
\cite{jimbo}.
In the natural basis of tensor products its matrix form reads
\be{2.2}
R=\left(
\arb{cccc} 1&0&0&0\\0&q^{-1}&0&0\\0&(1-q^{-2})& q^{-1}&0\\0&0&0&1 \are
\right).
\ee
If in addition one associates to each space $V_i$ a variable $x_i$, it is 
possible to define a 'spectral parameter' dependent R-matrix:
\be{2.3}
R(x):=\frac{qe^{x/2}R-q^{-1}e^{-x/2}PR^{-1}P}{qe^{x/2}-q^{-1}e^{-x/2}}=
\left(
\arb{cccc} 1&0&0&0\\0&b(x)&c_-(x)&0\\0&c_+(x)& b(x)&0\\0&0&0&1 \are
\right),
\ee 
where $P$ is the permutation operator in the sense of eqn.~(\ref{2.1})
\[
P_{ij}\, (v_i \ot v_j)= v_i \ot v_j, \quad v_{i,j}\in V_{i,j},
\]
and $x=x_i-x_j$. The Boltzmann weights read explicitly
\be{2.4}
b(x)=\frac{e^{x/2}-e^{-x/2}}{qe^{x/2}-q^{-1}e^{-x/2}},
\quad c_{\pm}(x)=\frac{e^{\pm x/2}(q-q^{-1})}{qe^{x/2}-q^{-1}e^{-x/2}}.
\ee
$R(x)$ satisfies the Yang Baxter equation
\be{2.5}
R_{12}(x_1-x_2)R_{13}(x_1-x_3)R_{23}(x_2-x_3)=
R_{23}(x_2-x_3)R_{13}(x_1-x_3)R_{12}(x_1-x_2).
\ee
One defines a monodromy matrix $T_0(\x,x_0)$ acting on the  
tensor product space $V=\bigotimes_{i=1}^N V_i$ and an additional 
auxiliary space $V_0 \simeq \mbox{${\bf C}^2$}$
\be{2.6}
T_0(\x,x_0):=R_{10}(x_1-x_0)R_{20}(x_2-x_0) \ldots R_{N0}(x_N-x_0).
\ee
However as we will see later it is more useful to work with the doubled 
monodromy matrix as proposed in \cite{skly} some years ago as an application 
of the 'reflection' equation introduced in \cite{cher}. We will use in the 
following the special type that has been introduced in \cite{kz1} and is 
given by
\be{2.7}
\tmod_0(\x,x_0):=R_{01}R_{02}\ldots R_{0N}
		R_{10}(x_1-x_0)R_{20}(x_2-x_0) \ldots R_{N0}(x_N-x_0).
\ee  
Its dependency on $V_0$ becomes obvious if $\tmod_0$ is written as a matrix 
w.r.t.~the auxiliary space:
\be{2.8}
\tmod_0=\left( \arb{cc} \amod & \bemod \\ \cmod & \dmod \are \right).
\ee
Equation (\ref{2.5}) implies the Yang-Baxter equation for $\tmod$
\be{2.9}
R_{ab}(v-u)\tmod_a(\x;u)R_{ba}\tmod_b(\x;v)=
\tmod_b(\x;v) R_{ba} \tmod_a(\x;u) R_{ba}(v-u)
\ee
giving the following commutation relations for the operators $\amod,\bemod$ 
and $\dmod$:
\ban{2.10}
\nonumber
\left[\bemod(\x;u),\bemod(\x;v)\right]&=&0,
\\
\nonumber
\amod(\x;u)\bemod(\x;v)&=& q\iv  b\iv (u-v) \bemod(\x;v)\amod(\x;u) 
\\
\nonumber
&&-\bemod(\x;u)\left[  q \iv  \frac{c_-(u-v)}{b(u-v)} \amod(\x;v) 
	+(1-q^{-2})\dmod(\x;v)\right],\\
\dmod(\x;u)\bemod(\x;v)&=&q b\iv (v-u)\bemod(\x;v)\amod(\x;u) -
	q \frac{c_-(v-u)}{b(v-u)}\bemod(\x;u)\dmod(\x;v),
\ean
(for a detailed proof see \cite{kz1}).

Analogous to the definition (\ref{2.7}) consider a further set of monodromy 
type matrices defined by
\be{2.11}
\tmod^Q(\x;i):=R_{01}R_{02}\ldots R_{0N}
		R_{10}(x_1-x_0)\ldots P_{i0} \ldots R_{N0}(x_N-(x_0+\kappa)),
\quad(i=1,\ldots,N).
\ee
where $\kappa$ is the arbitrary shift parameter having already appeared in 
eqn.~(\ref{1.1}).
They still have the block structure of (\ref{2.8}) but do no longer
depend on a parameter $x_0$ of the auxiliary space. The monodromy matrices
$\tmod$ and $\tmod^Q$ now fulfill another Yang-Baxter equation
\be{2.12}
R_{ba}(x_i-u)\,\tmod_b(\x';u)\,R_{ab}\,\tmod^Q_a(\x;i)=
\tmod^Q_a(\x;i)\, R_{ba}\, \tmod_b(\x;u)\, R_{ab}(x_i+\kappa-u).
\ee
Again we give some commutation rules relating their matrix elements:
\ba
\lefteqn{
\amod^Q(\x;i)\bemod(\x;u)= q\iv  b\iv (x_i+\kappa-u) \bemod(\x';u)
	\amod^Q(\x;i) 
}
\\
&& -\bemod^Q(\x;i)\left[  q \iv  \frac{c_-(x_i+\kappa-u)}{b(x_i+\kappa-u)} 
\amod(\x;u) 
	+(1-q^{-2})\dmod(\x;u)\right],
\ea
\be{2.13}
\dmod^Q(\x;i)\bemod(\x;u)=q b\iv (x_i-u)\bemod(\x';u)\dmod^Q(\x;i) -
	q \frac{c_-(u-x_i)}{b(u-x_i)}\bemod^Q(\x;i)\dmod(\x;u).
\ee
The first terms in eqn.~(\ref{2.10}) and (\ref{2.13}) are 
called 'wanted' resp. all others 'unwanted' ones. 
Now taking the Markov trace\footnote{It's asymmetric form results from the choice of normalization in 
eqn.~(\ref{2.2})} over $\tmod^Q$ gives the operator on the 
l.h.s.~of the difference equation (\ref{1.1})
\be{2.14}
Q(\x;i):=\tr_q\tmod^Q(\x;i)=\amod^Q(\x;i)+q^{-2}\dmod^Q(\x;i).
\ee
Denote by $\Omega$ the usual reference state ($\Omega=|1\rangle^{\ot N} $) 
and apply an arbitrary number $m$ of $\bemod$ operators thereto defining the  
following\\
{\bf Bethe ansatz vector:}
\be{2.16}
f(\x)=\sum_{\uu}B(\x;u_m) \ldots B(\x;u_1) \, \Omega \, g(\x;\uu),
\ee
where the summation 
over $\uu$ is specified by  
\be{2.17}
\sum_{\uu}=
\sum_{\stackrel{{\scriptstyle l_1 \in} \mbox{$\s \bf Z$}}
{u_1=\tilde{u}_1+l_1\kappa}}
\ldots
\sum_{\stackrel{{\scriptstyle l_m \in} \mbox{$\s \bf Z$}}
{u_m=\tilde{u}_m+l_m\kappa}}\mbox{($\tilde{\uu}$ arbitrary set of complex 
numbers)}
\ee
and the the function $g(\x;\uu)$ is defined by
\be{2.18}
g(\x;\uu)=\prod_{i,j}\psi(x_i-u_j)
\prod_{k<l}\tau(u_k-u_l).
\ee
{\bf Theorem:} The difference equation (\ref{1.1}) defined by eqn.\ 
({\ref{2.14}) is solved by the Bethe vectors (\ref{2.16}) if the functions 
$\psi(x)$ and $\tau(x)$ satisfy the difference equations:
\be{2.19}
q\iv b(x+\kappa) \, \psi(x+\kappa)=\psi(x),\quad q^2 \frac{\tau(x)}{b(x)}=
\frac{\tau(x-\kappa)}{b(-x+\kappa)}.
\ee
{\bf Remark:}   
As a variation of the solutions given
in \cite{resh} the following functions fulfill this conditions:
\be{2.20}
\psi(x)=
\frac{(q^2 e^x;e^{-\kappa})_\infty}
{(e^x;e^{-\kappa})_\infty},\quad
\tau(x)=
(1-e^x)\frac{(q^{-2} e^{x-\kappa};e^{-\kappa})_\infty}
{(q^{2}e^x;e^{-\kappa})_\infty},
\ee
where
\[
(z;p)_\infty:=\prod_{n=0}^{\infty}(1-zp^n).
\]
{\bf Proof:} 
We apply the operator $Q(\x;i)$ in its decomposition (\ref{2.14}) to $f(\x)$. 
Using the relations  
(\ref{2.11}) and (\ref{2.14}) one commutes  the operators $\amod^Q$ and 
$\dmod^Q$ to the right, where they act on the reference state $\Omega$ 
according to
\[ \amod^Q(\x;i)\Omega=\Omega, \quad \dmod^Q(\x,i)\Omega=0, \]
respectively
\[ 
\amod(\x;u)\Omega=\Omega, \quad \dmod(\x,u)\Omega=\prod_{j=1}^N 
	b(x_j-u) \, \Omega.
\]
The wanted term contribution of $\amod^Q$ reads
\ba
\lefteqn{
\amod(\x;i)\sum_{\uu}\bemod(\x;u_m)\ldots\bemod(\x;u_1)\,\Omega\, g(\x,\uu)
=}\\
&&= \sum_{\uu}\bemod(\x';u_m)\ldots\bemod(\x';u_1)\,\Omega\, 
\prod_{j=1}^{m}q \iv b\iv(x_i+\kappa-u_j) g(\x,\uu)=f(\x'),
\ea
where in the last step the quasi periodic property of $\psi$ 
(eqn.(\ref{2.19})) has been used. The $q^{-2}\dmod^Q$ wanted contribution 
vanishes due to the fact that $\dmod^Q(\x;i)\,\Omega=0$.

In a second step one has to verify that all other terms cancel each other 
under the sum (\ref{2.17}). Denote the unwanted 
terms obtained from $\amod^Q$ respectively $q^{-2}\dmod^Q$ that are 
proportional to $\bemod^Q(\x;i)\bemod(\x;u_{m-1})\ldots\bemod(\x;u_1)\,
\Omega$ by ${\rm uw}_{A,D}^{(i,j)}$. 
(They result when one commutes first 'unwanted' due to (\ref{2.13}) 
and then always wanted due to (\ref{2.10}).)

\lefteqn{{\rm uw}_A^{(i,m)}=}
\lefteqn{\left[-q\iv \frac{c_-(x_i+\kappa-u_m)}{b(x_i+\kappa-u_m)}
\prod_{k<m} q\iv b\iv(u_m-u_k)
-(1-q^{-2})\prod_{k<m}q b\iv(u_k-u_m)\prod_{j=1}^N q\iv b(x_j-u_m)\right]
}
\be{2.21}
\bemod^Q(\x;i)\bemod(\x;u_{m-1})\ldots\bemod(\x;u_1)
\,\Omega g(\x;\uu)
\ee
\lefteqn{{\rm uw}_D^{(i,m)}=}
\be{2.22}
 -q\iv\frac{c_-(u_m-x_i)}{b(u_m-x_i)}\prod_{k<m}q\,b\iv
	(u_k-u_m)\prod_{j=1}^N q\iv b(x_i-u_m)
\bemod^Q(\x;i)\bemod(\x;u_{m-1})\ldots\bemod(\x;u_1)\,\Omega g(\x;\uu).
\ee
Using the symmetry property 
$ \frac{c_-}{b}(-x)=- \frac{c_-}{b}(x)-(1-q^{-2})$ 
combine (\ref{2.22}) and the second term of (\ref{2.21}). 
Then both eqns. of (\ref{2.19}) are applied to this term and 
obviously this term cancels with the first one of (\ref{2.21}) 
under the sum (\ref{2.17}) which completes the proof.
%
%
\section{Bethe vectors and highest weight modules}
\setcounter{equation}{0}
The generators of $U_q[sl(2)]$ can be derived from the monodromy matrix 
$T_0(\x;u)$ (\ref{2.7}) in the limits $u \to \pm \infty$:
\ban{3.1}
\nonumber
T=\left(\arb{cc}T^{11}&T^{12}\\T^{21}&T^{22}\are\right):=
\lim_{u\to -\infty} T_0(\x;u)=
q^{-N}
&\left(\arb{cc} 1&0\\(q-q \iv)J_+&1\are\right)&
q^{\bf W}; 
\\
\tilde{T}:=
\lim_{u\to +\infty} T_0(\x;u)=
q^N q^{-\bf W}&\left(\arb{cc} 1&(q-q \iv)J_-\\0&1\are\right)&,
\ean
where ${\bf W}={\rm diag}\{W_1,W_2\}$ contains the Cartan elements. 
In order to prove the highest weight property of $f(\x)$ i.e.\ the 
statement
\be{3.2}
T^{21}f(\x)\propto J_+f(\x)=0,
\ee 
we introduce analogous to (\ref{3.1}) as 
a limit of $\tmod(\x;u)$ 
\be{3.3}
\tmod:=\tilde{T}^{-1}\,T.
\ee
First we show that $\tmod^{21}f(\x)=0$. The Yang-Baxter equation 
(\ref{2.12}) implies
\be{3.4}
\left[ \tmod^{21},\bemod(u)\right]=
(1-q^{-2})\left[\amod(u)\tmod^{22}-\tmod^{22}\dmod(u)\right].
\ee
Again due to the commutativity of the $\bemod$-operators it is sufficient to
consider the term proportional to $\bemod(u_m)\ldots\bemod(u_2)$. Because 
$\amod(u),\dmod(u)$ and $\tmod^{22}$ act diagonal on $\Omega$ it remains to 
show that
\[
\sum_{\uu}\left[\amod(u_1)-\dmod(u_1)\right] \, \Omega\, g(\x;\uu)=0,
\]
which follows directly from eqn.~(\ref{2.19}).  
Since $\tilde{T}\iv$ is an invertible operator eqn.\ ({3.3}) implies the 
statement (\ref{3.2}). 

The weights $\bf \omega$ of the Bethe vectors $f(\x)$
are defined by
\[
q^{\bf W} f(\x)=q^{\bf \omega}f(\x), \quad {\bf \omega}=(\omega_1,\omega_2).
\]
The commutation relations
\[
\amod\bemod(u)=q^{-2} \bemod(u)\amod;\quad \dmod\bemod(u)=q^2 \bemod(u)\dmod
\]
and eqn.~(\ref{3.3}) therefore imply
\[
T^{11}f(\x)=q^{N-m}f(\x)\quad \mbox{and} \quad T^{22}f(\x)=q^m f(\x),
\]
giving the weight vector ${\bf \omega}=(N-m,m)$ as expected.
\section{The higher ranked case $U_q[sl(n)]$}
\setcounter{equation}{0}
In this last section we briefly discuss the case of an $U_q[sl(n)]$ difference 
equation. (For a more detailed description of the nested Bethe ansatz method in
general we refer the reader to \cite{bkz1} and \cite{kz1})

Denote by $E_{ij}$ the unit matrices in $M_{n,n}(\bf C)$. The $U_q[sl(n)]$ 
R-matrix is then given by
\be{4.1}
R=\sum_{i} E_{ii} \ot E_{ii}+q\iv \sum_{i \neq j} E_{ii} \ot E_{jj}+
  (1-q^{-2}) \sum_{i>j} E_{ij} \ot E_{ji},
\ee
wheras the definitions for $R(x),\, T_0(\x;u)$ and $\tmod_0(\x;u)$ 
can directely be overtaken from the equations (\ref{2.3}),  (\ref{2.6}) and 
(\ref{2.7}). The latter two operators are considered now as $n \times n$ 
matrices; the commutation relations of their elements read in analogy to 
eqns.~(\ref{2.10}) and (\ref{2.13}).
\ban{4.2}
\nonumber
\amod(\x;u)\bemod_\gamma(\x;v) &=& q\iv b\iv (u-v) \bemod_\gamma(\x;v)
\amod(\x;u)\\
\nonumber
&&- q\iv \left[ \frac{c_-(u-v)}{b(u-v)} \bemod_\ga(\x;u) \amod(\x;v) +
(q-q\iv)\bemod_\alpha(\x;u)\dmod_{\alpha\ga}(\x;v) \right],\\
\nonumber
\dmod_{\beta\ga}(\x;u)\bemod_\de(\x;v) &=& q\, b\iv(v-u) [
\bemod_{\ga''}(\x;v) \dmod_{\beta'\de'}(\x;u)R^{\de'\ga'}_{\de\ga}(v-u)
R^{\ga''\beta}_{\beta'\ga'}\\
\nonumber &&
-c_-(v-u)R^{\ga'\beta}_{\beta'\ga} \bemod_{\ga'}(\x;u)\dmod_{\beta'\de}(\x;v)]
\\
\nonumber
\amod^Q(\x;i)\bemod_\gamma(\x;u) &=& q\iv b\iv (x_i+\kappa-u) 
\bemod_\gamma(\x';u)\amod^Q(\x;i)\\
\nonumber
&&- q\iv \left[ \frac{c_-(x_i+\kappa-u)}{b(x_i+\kappa-u)} \bemod^Q_\ga(\x;i) 
\amod(\x;u) +
(q-q\iv)\bemod^Q_\alpha(\x;i)\dmod_{\alpha\ga}(\x;u) \right],\\
\nonumber
\dmod^Q_{\beta\ga}(\x;i)\bemod_\de(\x;u) &=& q\, b\iv(u-x_i) [
\bemod_{\ga''}(\x';u) \dmod^Q_{\beta'\de'}(\x;i)R^{\de'\ga'}_{\de\ga}(u-x_i)
R^{\ga''\beta}_{\beta'\ga'}\\
\nonumber &&\left.
-c_-(u-x_i)R^{\ga'\beta}_{\beta'\ga} \bemod^Q_{\ga'}(\x;i)
\dmod_{\beta'\de}(\x;u)\right],
\ean
where the greek indices run from $2$ to $n$. The operators $Q(\x;i)$, which
define eqn.~({\ref{1.1}) are given by the $U_q[sl(n)]$ Markov trace
\be{4.3}
Q(\x;i):={\rm tr}_q \tmod^Q(\x;i)=\amod^Q(\x;i)+\sum_{\alpha=2}^n 
q^{-2(\alpha-1)}\dmod_{\alpha\alpha}^Q (\x;i),
\ee
Furtheron we denote the number of particles by $N_n$.
The Bethe vectors solving (\ref{1.1}) are created by the action of $N_{n-1}$ 
$\bemod$-Operators and read
\be{4.4}
f(\x)=\sum_{\uu} \bemod_{\beta_{N_{n-1}}}
(\x;u_{N_{n-1}}) \ldots \, 
\bemod_{\beta_1}(\x;u_1)\,
\Omega\, g^{\underline{\beta}}(\x,\uu),
\ee
where, in contrast to section 2, $g(\x;\uu)$ is a function with values in
$V^{(n-1)}=\ot^{N_{n-1}}{\bf C}^{(n-1)}$ given by the ansatz
\be{4.5}
g(\x;\uu)=\prod_{i,j}\psi(x_i-u_j) \prod_{k<l}\tau(u_k-u_l)\, f^{(n-1)}(\uu)
\ee
with functions $\psi(x)$ and $\tau(x)$ as given by (\ref{2.20}) and a 
(yet undetermined) function $f^{(n-1)}$ with values in $V^{(n-1)}$.
To prove eqn.~(\ref{1.1}) one has to apply $Q(\x;i)$ to $f(\x)$;
the 'wanted' contribution of $\amod^Q$ again produces the r.h.s.~of 
eqn.~(\ref{1.1}). On the other hand the 'unwanted' terms cancel exactly 
if $f^{(n-1)}$ satisfies the $n-1$ dimensional analogue of eqn.~(\ref{1.1}).
Therefore we repeat the ansatz (\ref{4.4}) for $f^{(n-1)}$ and all the
resulting subsequent Bethe ansatz levels, where consequently the number of 
$\bemod$-operators used at the $k$th level is denoted by $N_{n-k}$. Finally 
after $n-2$ steps the 
problem has been reduced to the $U_q[sl(2)]$ problem already solved in 
Section 2. 

The highest weight property of the Bethe vectors ({\ref{4.4}) 
is proved in a way parallel to Section 3. At some stages the higher ranked 
case is a little more involved, but those aspects have been already treated 
carefully in \cite{kz2}.

The resulting weight vector ${\bf \omega}$ then reads
\be{4.6}
\omega=(\omega_1,\ldots,\omega_n)=(N_n-N_{n-1},\ldots,N_2-N_1,N_1), 
\ee
again fulfilling the maximal weight condition
\be{4.7}
\omega_1 \geq \ldots \geq \omega_n \geq 0.
\ee
\section*{Summary}
Starting from the $U_q[sl(2)]$ R-matrix we derived a family of q-deformed 
discrete Knizhnik-Zamolodchikov equations and constructed solutions 
via the generalization of the algebraic Bethe ansatz as developed in 
\cite{bkz1}. These solutions have been shown to be of highest weight w.r.t. 
the underlying quantum group structure. 
Using the variant of the nested Bethe ansatz method we extended the results 
to the higher ranked symmetry of $U_q[sl(n)]$.
An application of these results can be found in \cite{bfk}.
%
\paragraph{Acknoledgements.} 
The author would like to thank H.~Babujian and 
M.~Karowski for numerous helpful and stimulating discussions.
\pagebreak 


\begin{thebibliography}{99}
\bibitem{bkz1} 
H.~Babujian, M.~Karowski and A.~Zapletal, 
J.~Phys.~A Math.~Gen.~{\bf 30} (1997) 6425.
%
\bibitem{tv} V.~Tarasov and A.~Varchenko, Invent.~math.~{\bf 128} (1997) 501.
\bibitem{FR} 
I.~Frenkel and N.Yu.~Reshetikhin, Commun.~Math.~Phys.~{\bf 146} (1992) 1.
%
\bibitem{resh} 
N.Yu.~Reshetikhin, Lett.~Math.~Phys.~{\bf26} (1992) 153.
%
\bibitem{nv} 
L.~Mezincescu and R.I.~Nepomechie, Mod.~Phys.~Lett.~{\bf A6} (1991)  
2497;\\
C.~Destri and H.J.~de Vega, Nucl.~Phys.~{\bf B361}
(1992) 361.
%
\bibitem{kz1} 
M.~Karowski and A.~Zapletal, Nucl.~Phys.~{\bf B419} (1994) 567.
%
\bibitem{jimbo}
M.~Jimbo, Lett.~Math.~Phys.~{\bf 10} (1985) 63.
\bibitem{skly} E.K.~Sklyanin, J.~Phys.~{\bf A21} (1988) 2375.
%
\bibitem{cher} I.~Cherednik, Theor.~Mat.~Fiz.~{\bf 61} (1984) 35.
%
\bibitem{Be} 
H.~Bethe, Z.~Phys.~{\bf 71} (1931) 205.
%
\bibitem{TF} 
L.A.~Takhtadzhan and L.D.~Faddeev, Russian Math.~Surveys {\bf 34} (1979) 186.
\bibitem
{kz2} M.~Karowski and A.~Zapletal, J.~Phys.~A Math.~Gen. {\bf 27} (1994) 7419.
\bibitem
{bfk} H.~Babujian, A.~Fring, M.~Karowski and A.~Zapletal, Exact Form Factors in Integrable Field Theories: the Sine-Gordon Model, hep-th/9805185.
\end{thebibliography}
\end{document}